\documentclass[12pt,a4paper]{article}
\usepackage{amssymb, amsmath}  
\begin{document}
\title{\textbf{Stability of the LCD Model}}

\author{Zhenting Hou$^{1}$,
 Li Tan$^{1}${\thanks{Correspondence to: Li Tan; E-mail:tltanli@yahoo.com.cn}, Dinghua Shi$^{1,2}$}
 \\  $^{1}$School of Mathematics, Central South University, \\Changsha, Hunan, 410075, China\\
   $^{2}$Department of Mathematics, Shanghai University,\\ Shanghai, 200444, China}

\maketitle

{\textbf{Abstract:}} In this paper, first-passage probability of
Markov chains is used to get a strict proof of the existence of
degree distribution of the LCD model presented by Bollob\'{a}s
(Random Structures and Algorithms 18(2001)). Also, a precise
expression of degree distribution is presented.


\newtheorem{theorem}{Theorem}
\newtheorem{remark}{Remark}
\newtheorem{defi}[theorem]{Definition}
\newtheorem{Lemma}{Lemma}
\bigskip

{\it Introduction}. Networks exist in every expects of nature and
society. Most of the systems, e.g. WWW, movie actor collaboration
network, transportation network and so on, can be described as
complex networks. Thus, complex networks becomes one of the
important research areas in recent years.

Two important discoveries about complex networks is the small-world
property[11] and the scale-free property, proposed by Barab\'{a}si
and Albert[1]. They had a research on World-Wide Web and found that
the degree distribution of the WWW is not Poisson, as random
networks and small-world networks shows, but a power-law
distribution. In[1], they presented the approximate expression of
network degree and used computer experiments to get $\gamma=2.9\pm
0.1$. It is stated that after many steps the proportion of vertices
with degree $d$ should obey a power law $P(d)\propto d^{-\gamma}$.
In[4] they stated an evolving network(BA model) with growth and
preferential attachment. They used the mean-field method to compute
for the degree distribution, and deduced $\gamma=3$ heuristically.
They thought that the probability of one vertex to get edges equals
$m\Pi(k_i)$. That is to say, it is in direct proportion to degree of
the vertex and the number of edges added at a unit time
simultaneously. They also pointed out that most real networks are
scale-free[2].

In order to make a better understanding of the phenomenon in real
networks, many scholars did researches on BA model, they also had
extensions to BA model and got a lot of useful results.

Krapivsky[12] let $N_k(t)$ be the number of vertices in network with
degree $k$ at time $t$. They deduced the rate equations:
$\frac{dN_k(t)}{dt}=m\frac{(k-1)N_{k-1}(t)-kN_k(t)}{\sum\limits_{k}kN_k(t)}+\delta_{km}$
and used $\lim\limits_{t\rightarrow\infty}E[N_k(t)/t]=P(k)$ to show
that $P(k)=\frac{4}{k(k+1)(k+2)}$ for the BA model with $m=1$.

Dorogovtsev[8] considered  $k_i(t)$, the degree of the vertex added
at time $i$ evolved at time $t$. They used $P(k,i,t)$ to express the
probability of vertex $i$ having $k$ edges at time $t$. And let the
average degree $P(k,t)=\frac{1}{t}\sum\limits_{i=1}^{t}P(k,i,t)$ to
be the network degree at time $t$. They presented the attraction
model, where each new vertex has an initial attraction degree $a$.
Thus, $k=a+q$ is the total degree with $q$ being the in-degree of
the vertex. They got the master equation:
$P(k,i,t+1)=\frac{k-1}{2t}P(k-1,i,t)+(1-\frac{k}{2t})P(k,i,t)$ and
the stationary distribution: $P(k)=\frac{2m(m+1)}{k(k+1)(k+2)}$ with
$a=m$.

These two definitions of network degree is equivalent, for we have $
E[N_k(t)]=\sum\limits_{i=1}^{t}P(k,i,t).$

As to BA model, Bollob\'{a}s[6] once raised the following questions:
{\it The first is getting started: how do we take probabilities
proportional to the degrees when these are all zero? The second
problem is with the preferential attachment rule itself. For the BA
model, or for any other BA-like model, how to find a rigorous
preferential attachment scheme for $m\geq 2$?}

Aimed at this questions, Bollob\'{a}s[5] presented a modified BA
model, called LCD model.

As Bollob\'{a}s[5] stated:{\it We start with the case $m=1$.
Consider a fixed sequence of vertices $v_1,v_2,\cdots$. We shall
inductively define a random graph process $(G_1^t)_{t\geq 0}$ so
that $G_1^t$ is a directed graph on $\{v_i:1\leq i\leq t\}$, as
follows. Start with $G_1^0$ the graph with no vertices, or with
$G_1^1$ the graph with one vertex and one loop. Given $G_1^{t-1}$,
form $G_1^{t}$ by adding the vertex $v_t$ together with a single
edge directed from $v_t$ to $v_i$, where $i$ is chosen randomly with
\begin{equation}
P\{i=s\}=\left\{
\begin{array}{cc}
\frac{dG_1^{t-1}(v_s)}{2t-1}&1\leq s\leq t-1,\\
\frac{1}{2t-1}&s=t.
\end{array}
\right. \label{Pis}
\end{equation}

For $m>1$ add $m$ edges from $v_t$ one at a time. We define the
process $(G_m^t)_{t\geq 0}$ by running the process $(G_1^{t})$ on a
sequence $v'_1,v'_2,\cdots$; the graph $G_m^{t}$ is formed from
$G_1^{mt}$ by identifying the vertices $v'_1,v'_2,\cdots,v'_m$ to
form $v_1$, identifying $v'_{m+1},v'_{m+2},\cdots,v'_{2m}$ to form
$v_2$, and so on.}

Bollob\'{a}s[5] found that, although the graph process $\{G_1^t\}$
is dynamic, there is a static description (Linearized Chord Diagram)
of the distribution of the graph at a particular time. In[7], they
studied $\sharp_m^n(d)$, the number of vertices with indegree $d$
and strictly proved the following equation
\begin{equation}
E(\sharp_m^n(d))\sim\frac{2m(m+1)n}{(d+m)(d+m+1)(d+m+2)}.
\label{Emn}
\end{equation}
Then, with Azuma-Hoeffding inequality[10], they got the following
theorem.

{\it Theorem  Let $m\geq 1$ be fixed, and consider the random graph
process $(G_m^n)_{n\geq 0}$, writing $\sharp_m^n(d)$ for the number
of vertices of $G_m^n$ with indegree equal to $d$, i.e., with total
degree $m+d$. Let
\begin{eqnarray*}
\alpha_{m,d}=\frac{2m(m+1)}{(d+m)(d+m+1)(d+m+2)}
\end{eqnarray*}
and let $\epsilon$ be fixed. Then with probability tending to 1 as
$n\rightarrow\infty$we have
\begin{eqnarray*}
(1-\epsilon)\alpha_{m,d}\leq\frac{\sharp_m^n(d)}{n}\leq(1+\epsilon)\alpha_{m,d}
\end{eqnarray*}
for every d in the range $0\leq d \leq n^\frac{1}{15}.$}

In fact, (\ref{Emn}) already proved that $\sharp_m^n(d)/n$
convergent to $\alpha_{m,d}$ with first order. The theorem is used
to prove that $\sharp_m^n(d)/n$ tends to $\alpha_{m,d}$ with
probability 1. If one is only interested in the degree distribution
of the LCD model, (\ref{Emn}) is enough. Also, Bollob\'{a}s[7]
changed the BA model and the method used in his paper is applicable
only to particular models and it is difficult to have extensions.

In this paper, we apply first-passage probability[9] to LCD model.
We used the Markov theory and got a rigorous proof the existence of
the steady-state degree distribution of the LCD model. And finally
presented the precise expression of the degree distribution.

{\it Stability}. Instead of studying martingale
$X_I=E[\sharp_m^T(d)|G_m^I] (where 0\leq I\leq T)$, we consider
$k_I(T)$, the degree of vertex $I$ at time $T$. Following
Dorogovtsev[8], $k_I(T)$ is a random variable, and let
$P(k,I,T)=P\{k_I(T)=k\}$ the probability of vertex $I$ having $k$
edges at time $T$, also take the average degree
$P(k,T)=\frac{1}{T}\sum\limits_{I=1}^{T}P(k,I,T)$ to be the
definition of the network degree at time $T$. For fixed $T$,
$k_I(T)$ is a random variable and for variable $T$, $k_I(T)$ is a
nonhomogeneous Markov chain[13]. The adding of vertex $I$ can be
divided into $m$ steps, from $(I-1)m+1$ to $(I-1)m+m$. The
state-transition probability of this Markov chain is given by:
\begin{equation}
\begin{array}{lll}
& &P\{k_{I}(T+1)=l|k_{I}(T)=k\}\\
&=&\left\{
\begin{array}{cc}
P_{k,j}^{T+1}&l=k+j,j=0,1,\cdots,m\\
0&otherwise.
\end{array}
\right.
\end{array}
\label{PIT}
\end{equation}
where $k=m,m+1,\cdots,(T-I+2)m$ and $P_{k,j}^{T+1}$ is the
probability of vertex $I$ with degree $k$ at $T$ increase by $j$ at
$T+1$. With (\ref{Pis}), it follows that
\begin{eqnarray*}
P_{k,0}^{T+1}&=&\prod\limits_{r=1}^{m}(1-\frac{k}{2Tm+2(r-1)+1})\\
P_{k,1}^{T+1}&=&\sum\limits_{b=1}^{m}[\prod\limits_{q=1}^{b-1}(1-\frac{k}{2Tm+2(q-1)+1})(\frac{k}{2Tm+2(b-1)+1})\\
         & &\prod\limits_{r=b+1}^{m}(1-\frac{k+1}{2Tm+2(r-1)+1})]\\
&\vdots&\\
P_{k,m}^{T+1}&=&\prod\limits_{r=1}^{m}\frac{k+r-1}{2Tm+2(r-1)+1}
\end{eqnarray*}

Consider the first-passage probability $f(k,I,S)$ of the Markov
chain $k_I(T)$. $f(k,I,S)=P\{k_I(S)=k, k_I(L)\neq k,
L=1,2,\cdots,S-1\}$. Based on the concept and techniques of Markov
chain, the relationship between the first-passage probability and
the probability of vertex degrees are as follows:

\begin{Lemma}
When $1\leq I<S, k>m$, we have
\begin{equation}
f(k,I,S)=\sum\limits_{j=1}^{m}P(k-j,I,S-1)P_{k-j,j}^{S},
\label{fkIS}
\end{equation}

When $1\leq I<T, k>m$, we have
\begin{equation}
\begin{array}{lll}
P(k,I,T)&=&\sum\limits_{S=I}^{T}f(k,I,S)
\prod\limits_{A=S+1}^{T}\prod\limits_{q=1}^{m}(1-\frac{k}{2(A-1)m+2(q-1)+1}).
\end{array}
\label{PkIT}
\end{equation}
\end{Lemma}
{\it Proof:} From the construction of LCD model, the degree of
vertex $I$ is always nondecreasing. At $S$, the degree of vertex $I$
can increase by $1, 2, \cdots,$ or $m$. With the Markov property, we
have
\begin{eqnarray*}
f(k,I,S)&=&P\{k_I(S)=k, k_I(L)\neq k, L=1,2,\cdots,S-1\}\\
        &=&\sum\limits_{j=1}^{m}P\{k_I(S)=k, k_I(S-1)=k-j\}P_{k-j,j}^{S}\\
        &=&\sum\limits_{j=1}^{m}P(k-j,I,S-1)P_{k-j,j}^{S}.
\end{eqnarray*}
Thus (\ref{fkIS}) is established.

Second, observe that the earliest time for the degree of vertex $I$
to reach $k$ is at step $I$, and the latest time to do so is at step
$T$. After this vertex degree becomes $k$, it will not increase any
more. \\This completes the proof. $\square$

\begin{Lemma}
When $1<I=S, k=m$, we have
\begin{equation}
f(m,I,I)=\prod\limits_{r=1}^{m}(1-\frac{r}{2(I-1)m+2(r-1)+1}).
\label{fmII}
\end{equation}

When $1<I=S, m<k<2m$, we have
\begin{equation}
\begin{array}{lll}
f(k,I,I)&=&\sum\limits_{r=k-m}^{m}\hat{P}(k-(m-r)-2,I,(I-1)m+r-1)\frac{k-(m-r)-1}{2(I-1)m+2(r-1)+1}\\
         &&\prod\limits_{q=r+1}^{m}(1-\frac{k-(m-r)+(q-r)}{2(I-1)m+2(q-1)+1}),
\end{array}
\label{fkII}
\end{equation}
Where $\hat{P}(k,I,(I-1)m+r)$ means the probability of the degree of
vertex $I$ to be $k$ at the $r$th step of adding in.

When $1<I=S, k=2m$, we have
\begin{equation}
f(2m,I,I)=\prod\limits_{r=1}^{m}(\frac{2(r-1)+1}{2(I-1)m+2(r-1)+1}).
\label{f2mII}
\end{equation}
\end{Lemma}

{\it Proof:} For $1<I=S, k=m$, it means vertex $I$ have no loops
after the first $m$ steps. For $1<I=S, k=2m$, it means vertex $I$
have $m$ loops after the first $m$ steps. Thus (\ref{fmII}) and
(\ref{f2mII}) are obtained.

In order to prove (\ref{fkII}), note that from step $(I-1)m+1$ to
$(I-1)m+m$ are all the adding of vertex $I$. At these $m$ steps, the
degree of vertex $I$ will at least increase by 1. In order for the
degree of vertex $I$ first reach $k$ at $(I-1)m+m$, then at
$(I-1)m+r$ the degree of vertex $I$ is at most $k-(m-r)$. Thus
(\ref{fkII}) is obtained.\\
This completes the proof. $\square$

{\it Stolz theorem: In sequence $\{\frac{x_n}{y_n}\}$, assume that
$\{y_n\}$ is monotone increasing sequence with
$y_n\rightarrow\infty$, if
$\lim\limits_{n\rightarrow\infty}\frac{x_{n+1}-x_n}{y_{n+1}-y_n}=l$
exists, where $-\infty\leq l \leq\infty$, then
$\lim\limits_{n\rightarrow\infty}\frac{x_n}{y_n}=l$.}\\
{\it Proof:} See[14].

\begin{Lemma}
$\lim\limits_{T\rightarrow\infty}P(m,T)$ exists, and
\begin{equation}
P(m)\triangleq\lim\limits_{T\rightarrow\infty}P(m,T)=\frac{2}{m+2}>0.
\label{Pm}
\end{equation}
\end{Lemma}
{\it Proof:} From the construction of LCD model, it follows that
\begin{eqnarray*}
&&P(m,I,T)=f(m,I,I)\prod\limits_{A=I+1}^{T}\prod\limits_{q=1}^{m}(1-\frac{m}{2(A-1)m+2(q-1)+1})\\
&=&\prod\limits_{r=1}^{m}(1-\frac{r}{2(I-1)m+2(r-1)+1})
         \prod\limits_{A=I+1}^{T}\prod\limits_{q=1}^{m}(1-\frac{m}{2(A-1)m+2(q-1)+1}).
\end{eqnarray*}
In particular, $P(m,1,T)=0$. So,
\begin{eqnarray*}
&&P(m,T)=\frac{1}{T}\sum\limits_{I=1}^{T}P(m,I,T)=\frac{1}{T}\sum\limits_{I=2}^{T}P(m,I,T)\\
      &=&\frac{1}{T}\sum\limits_{I=2}^{T}\prod\limits_{r=1}^{m}(1-\frac{r}{2(I-1)m+2(r-1)+1})
         \prod\limits_{A=I+1}^{T}\prod\limits_{q=1}^{m}(1-\frac{m}{2(A-1)m+2(q-1)+1})\\
      &=&\frac{1}{T}\prod\limits_{A=3}^{T}\prod\limits_{q=1}^{m}
         (1-\frac{m}{2(A-1)m+2(q-1)+1})\{\prod\limits_{r=1}^{m}(1-\frac{r}{2m+2(r-1)+1})\\
      & &+\sum\limits_{I=3}^{T}\prod\limits_{r=1}^{m}(1-\frac{r}{2(I-1)m+2(r-1)+1})
         \prod\limits_{B=3}^{I}\prod\limits_{h=1}^{m}(1-\frac{m}{2(B-1)m+2(h-1)+1})^{-1}\}
\end{eqnarray*}
Next, let
\begin{eqnarray*}
\begin{array}{lll}
X_N&=&\prod\limits_{r=1}^{m}(1-\frac{r}{2m+2(r-1)+1})+\sum\limits_{I=3}^{N}\prod\limits_{r=1}^{m}(1-\frac{r}{2(I-1)m+2(r-1)+1})
      \prod\limits_{B=3}^{I}\prod\limits_{h=1}^{m}(1-\frac{m}{2(B-1)m+2(h-1)+1})^{-1},
\end{array}
\end{eqnarray*}
and
\begin{eqnarray*}
Y_N=N\prod\limits_{A=3}^{N}\prod\limits_{q=1}^{m}(1-\frac{m}{2(A-1)m+2(q-1)+1})^{-1}>0.
\end{eqnarray*}
Thus, it follows that
\begin{eqnarray*}
X_{N+1}-X_N=\prod\limits_{r=1}^{m}(1-\frac{r}{2Nm+2(r-1)+1})
            \prod\limits_{B=3}^{N+1}\prod\limits_{h=1}^{m}(1-\frac{m}{2(B-1)m+2(h-1)+1})^{-1},
\end{eqnarray*}

\begin{eqnarray*}
Y_{N+1}-Y_N=(N+1-N\prod\limits_{r=1}^{m}(1-\frac{m}{2Nm+2(r-1)+1}))\\
            \prod\limits_{A=3}^{N+1}\prod\limits_{q=1}^{m}(1-\frac{m}{2(A-1)m+2(q-1)+1})^{-1}>0.
\end{eqnarray*}
Since $Y_N>0$ and $Y_{N+1}-Y_N>0$, $\{Y_N\}$ is a strictly monotone
increasing nonnegative sequence, hence $Y_N\rightarrow\infty$. Also,
by assumption,
\begin{eqnarray*}
\frac{X_{N+1}-X_N}{Y_{N+1}-Y_N}&=&\frac{\prod\limits_{r=1}^{m}(2Nm+2(r-1)+1)+o(N^m)}
   {\prod\limits_{r=1}^{m}(2Nm+2(r-1)+1)+mN\sum\limits_{i=1}^{m}\prod\limits_{r\neq i}^{m}
   (2Nm+2(r-1)+1)+o(N^m)}\\
   &\rightarrow&\frac{2}{m+2},(N\rightarrow\infty).
\end{eqnarray*}
With Stolz theorem, we have
\begin{eqnarray*}
P(m)\triangleq\lim\limits_{T\rightarrow\infty}P(m,T)=\lim\limits_{N\rightarrow\infty}\frac{X_N}{Y_N}
  =\lim\limits_{N\rightarrow\infty}\frac{X_{N+1}-X_N}{Y_{N+1}-Y_N}
  =\frac{2}{m+2}>0.
\end{eqnarray*}
This completes the proof.$\square$

\begin{Lemma}
When $k>m$, if $\lim\limits_{T\rightarrow\infty}P(k-1,T)$ exists,
then $\lim\limits_{T\rightarrow\infty}P(k,T)$ also exists and
\begin{equation}
P(k)\triangleq\lim\limits_{T\rightarrow\infty}P(k,T)=\frac{k-1}{k+2}P(k-1)>0
\label{Pk}
\end{equation}
\end{Lemma}

{\it Proof:} With (\ref{fkIS}) and (\ref{PkIT}) of Lemma 1, we have
\begin{eqnarray*}
P(k,T)&=&\frac{1}{T}\sum\limits_{I=1}^{T-1}P(k,I,T)+\frac{1}{T}P(k,T,T)\\
&=&\frac{1}{T}\sum\limits_{I=1}^{T-1}\sum\limits_{S=I}^{T}f(k,I,S)
\prod\limits_{A=S+1}^{T}\prod\limits_{q=1}^{m}(1-\frac{k}{2(A-1)m+2(q-1)+1})+\frac{1}{T}P(k,T,T)\\
&=&\frac{1}{T}\sum\limits_{I=1}^{T-1}f(k,I,I)
\prod\limits_{A=I+1}^{T}\prod\limits_{q=1}^{m}(1-\frac{k}{2(A-1)m+2(q-1)+1})+\frac{1}{T}P(k,T,T)\\
&&+\frac{1}{T}\sum\limits_{I=1}^{T-1}\sum\limits_{S=I+1}^{T}\sum\limits_{j=1}^{m}P(k-j,I,S-1)P_{k-j,j}^{S}
\prod\limits_{A=S+1}^{T}\prod\limits_{q=1}^{m}(1-\frac{k}{2(A-1)m+2(q-1)+1})\\
&=&f(T)+\frac{1}{T}\sum\limits_{I=1}^{T-1}\sum\limits_{S=I+1}^{T}\sum\limits_{b=1}^{m}\frac{k-1}{2(S-1)m+2(b-1)+1}P(k-1,I,S-1)\\
&&\prod\limits_{A=S+1}^{T}\prod\limits_{q=1}^{m}(1-\frac{k}{2(A-1)m+2(q-1)+1})+\frac{1}{T}P(k,T,T)+O(\frac{1}{T})\\
\end{eqnarray*}

Consider $f(T)$ firstly, obviously, $f(k,1,1)=\delta_{k,2m}$. Also,
when $1\leq I=S, k>2m$, then $f(k,I,I)=0$. With (\ref{fkII}) of
Lemma 2, when $1<I=S, m<k<2m$, we have
\begin{eqnarray*}
f(T)&=&\frac{1}{T}\sum\limits_{I=1}^{T-1}f(k,I,I)
\prod\limits_{A=I+1}^{T}\prod\limits_{q=1}^{m}(1-\frac{k}{2(A-1)m+2(q-1)+1})\\
&<&\frac{1}{T}\sum\limits_{I=2}^{T-1}f(k,I,I)<\frac{1}{T}\sum\limits_{I=2}^{T-1}\sum\limits_{r=k-m}^{m}
   \frac{k-(m-r)-1}{2(I-1)m+2(r-1)+1}\\
&<&\frac{1}{T}\sum\limits_{I=2}^{T-1}\frac{k-1}{2(I-1)}<(k-1)\frac{\ln
T+\gamma}{T}.
\end{eqnarray*}
Take limits on both sides, then
$\lim\limits_{T\rightarrow\infty}f(T)=0$ with $m<k<2m$.

Then, with (\ref{f2mII}) of Lemma 2, when $1<I=S, k=2m$, we have
\begin{eqnarray*}
f(T)<\frac{1}{T}(1+\sum\limits_{I=2}^{T-1}f(2m,I,I))
   <\frac{1}{T}(1+\sum\limits_{I=2}^{T-1}\frac{1}{(I-1)^m}).
\end{eqnarray*}
Take limits on both sides, then
$\lim\limits_{T\rightarrow\infty}f(T)=0$ with $k=2m$.

Thus, when $k>m$, we have $\lim\limits_{T\rightarrow\infty}f(T)=0$.

Now, consider the second item.
\begin{eqnarray*}
&&\frac{1}{T}\sum\limits_{I=1}^{T-1}\sum\limits_{S=I+1}^{T}\sum\limits_{b=1}^{m}\frac{k-1}{2(S-1)m+2(b-1)+1}P(k-1,I,S-1)\\
&&\prod\limits_{A=S+1}^{T}\prod\limits_{q=1}^{m}(1-\frac{k}{2(A-1)m+2(q-1)+1})\\
&=&\frac{1}{T}\sum\limits_{S=2}^{T}\sum\limits_{b=1}^{m}\frac{k-1}{2(S-1)m+2(b-1)+1}\sum\limits_{I=1}^{S-1}P(k-1,I,S-1)\\
&&\prod\limits_{A=S+1}^{T}\prod\limits_{q=1}^{m}(1-\frac{k}{2(A-1)m+2(q-1)+1})\\
&=&\frac{1}{T}\sum\limits_{S=2}^{T}\sum\limits_{b=1}^{m}\frac{(k-1)(S-1)}{2(S-1)m+2(b-1)+1}P(k-1,S-1)
\prod\limits_{A=S+1}^{T}\prod\limits_{q=1}^{m}(1-\frac{k}{2(A-1)m+2(q-1)+1})\\
&=&\frac{1}{T}\prod\limits_{A=3}^{T}\prod\limits_{q=1}^{m}(1-\frac{k}{2(A-1)m+2(q-1)+1})
\{\sum\limits_{b=1}^{m}\frac{k-1}{2m+2(b-1)+1}P(k-1,1)\\
&&+\sum\limits_{S=3}^{T}\sum\limits_{b=1}^{m}\frac{(k-1)(S-1)}{2(S-1)m+2(b-1)+1}P(k-1,S-1)
\prod\limits_{B=3}^{S}\prod\limits_{h=1}^{m}(1-\frac{k}{2(B-1)m+2(h-1)+1})^{-1}\}\\
\end{eqnarray*}
Next, let
\begin{eqnarray*}
X_N&=&\sum\limits_{b=1}^{m}\frac{k-1}{2m+2(b-1)+1}P(k-1,1)
+\sum\limits_{S=3}^{N}\sum\limits_{b=1}^{m}\frac{(k-1)(S-1)}{2(S-1)m+2(b-1)+1}P(k-1,S-1)\\
&&\prod\limits_{B=3}^{S}\prod\limits_{h=1}^{m}(1-\frac{k}{2(B-1)m+2(h-1)+1})^{-1}\},
\end{eqnarray*}
and
\begin{eqnarray*}
Y_N=N\prod\limits_{A=3}^{N}\prod\limits_{q=1}^{m}(1-\frac{k}{2(A-1)m+2(q-1)+1})^{-1}>0.
\end{eqnarray*}
Thus, it follows that
\begin{eqnarray*}
X_{N+1}-X_N&=&\sum\limits_{b=1}^{m}\frac{(k-1)N}{2Nm+2(b-1)+1}P(k-1,N)
\prod\limits_{B=3}^{N+1}\prod\limits_{h=1}^{m}(1-\frac{k}{2(B-1)m+2(h-1)+1})^{-1},
\end{eqnarray*}

\begin{eqnarray*}
Y_{N+1}-Y_N&=&\{N+1-N\prod\limits_{r=1}^{m}(1-\frac{k}{2Nm+2(r-1)+1})\}\\
 &&\prod\limits_{A=3}^{N+1}\prod\limits_{q=1}^{m}(1-\frac{k}{2(A-1)m+2(q-1)+1})^{-1}>0.
\end{eqnarray*}

Since $Y_N>0$ and $Y_{N+1}-Y_N>0$, $\{Y_N\}$ is a strictly monotone
increasing nonnegative sequence, hence $Y_N\rightarrow\infty$.
\begin{eqnarray*}
\frac{X_{N+1}-X_N}{Y_{N+1}-Y_N}&=&\frac{\sum\limits_{b=1}^{m}\frac{(k-1)N}{2Nm+2(b-1)+1}P(k-1,N)}
{N+1-N\prod\limits_{r=1}^{m}(1-\frac{k}{2Nm+2(r-1)+1})}\rightarrow\frac{k-1}{k+2}P(k-1),
(N\rightarrow\infty).
\end{eqnarray*}
With Stolz theorem, we have
\begin{eqnarray*}
P(k)&=&\lim\limits_{T\rightarrow\infty}(f(T)+\frac{1}{T}P(k,T,T)+O(\frac{1}{T}))+\lim\limits_{N\rightarrow\infty}\frac{X_N}{Y_N}\\
&=&\lim\limits_{N\rightarrow\infty}\frac{X_N}{Y_N}
=\lim\limits_{N\rightarrow\infty}\frac{X_{N+1}-X_N}{Y_{N+1}-Y_N}
=\frac{k-1}{k+2}P(k-1)>0.
\end{eqnarray*}
This completes the proof. $\square$

\begin{theorem}
The steady-state degree distribution of the LCD model exists, and is
given by
\begin{equation}
P(k)=\frac{2m(m+1)}{k(k+1)(k+2)}\thicksim 2m^2k^{-3}>0.
\end{equation}
\end{theorem}
{\it Proof:} By mathematical induction, if follows from Lemmas 3 and
4 that the steady-state degree distribution of the LCD model exists.
Then, solving (\ref{Pk}) iteratively, we obtain
\begin{eqnarray*}
P(k)=\frac{k-1}{k+2}P(k-1)=\frac{k-1}{k+2}\frac{k-2}{k+1}\frac{k-3}{k}P(k-3).
\end{eqnarray*}
By continuing the process till $k=3+m$, we get
\begin{eqnarray*}
P(k)=\frac{2m(m+1)}{k(k+1)(k+2)}\thicksim 2m^2k^{-3}>0.
\end{eqnarray*}
This completes the proof.$\square$

{\it Conclusion} The method of first-passage probability developed
by Hou et al.[9] has great generality. We apply this method to the
LCD model, and got the precise expression of degree distribution. Of
course, allowing loops and multiple edges made the LCD model
different from the BA model. Thus, the transition probability of
Markov chains changed a lot. And the method is quite different from
[9]. Moreover, the method presented in this paper can be extended to
the attraction model[8] and other models with multiple edges.

\end{document}